\documentclass[12pt,a4paper]{amsart}
\usepackage{amsmath,amssymb,amsfonts,amsthm,enumerate}

\newtheorem{lemma}{Lemma}

\theoremstyle{remark}
\theoremstyle{remark}\newtheorem{remark}{Remark}
\theoremstyle{plain}\newtheorem{theorem}{Theorem}
\theoremstyle{plain}\newtheorem*{thmcite}{Theorem A}
\theoremstyle{plain}\newtheorem*{thmcite1}{Theorem B}

\newcommand{\R}{\mathbb{R}}

\newcommand{\cs}{^*}

\newcommand{\biz}{\textbf{Proof.\ }}
\newcommand{\ri}{\right)}
\newcommand{\lef}{\left(}

\newcommand{\lk}{\left\lbrace}
\newcommand{\rk}{\right\rbrace}

\newcommand{\lsz}{\left[}
\newcommand{\rsz}{\right]}

\newcommand{\tv}{\rightarrow \infty}

\newcommand{\al}{\alpha}

\newcommand{\veps}{\varepsilon}

\title[Asymptotics of renewal-like equations]{Asymptotics of a renewal-like recursion and an integral equation}
 
\author{\'Agnes Backhausz}
\address{Department of Probability Theory and Statistics\\
 E\"otv\"os Lor\'and University\\P\'azm\'any P.~s.\ 1/C, H-1117
 Budapest, Hungary} 
\email{agnes@cs.elte.hu}
\author{Tam\'as F.~M\'ori}
\address{Department of Probability Theory and Statistics\\
E\"otv\"os Lor\'and University\\P\'azm\'any P.~s. 1/C, 
H-1117 Budapest, Hungary}
\email{moritamas@ludens.elte.hu}
\dedicatory{\upshape Department of Probability Theory and Statistics,\\
E\"otv\"os Lor\'and University\\ P\'azm\'any P.~s. 1/C,
H-1117 Budapest, Hungary\\                         
\textit{E-mail address:} 
\texttt{agnes@cs.elte.hu, moritamas@ludens.elte.hu} }

\thanks{The European Union and the European Social Fund have provided
financial support to the project under the grant agreement no. T\'AMOP
4.2.1./B-09/KMR-2010-0003.}
\date{21 May 2012}
\begin{document}

\begin{abstract}
We consider a renewal-like recursion and prove that the solution 
is polynomially decaying under suitable conditions. We prove similar
results for the corresponding integral equation. In both cases
coefficients and functions are of more general form than in the
classic cases. 
\end{abstract}
\maketitle
\thispagestyle{empty}
\section{Introduction}
In this paper we examine the asympotics of a renewal-like recursion
and a similar integral equation. The motivation comes from probability
theory; more precisely, in a random model of publication activity
\cite{publ} the asymptotic distribution of the weights of the
authors satisfy such equations.  

The recursion is of the form
\begin{equation}
x_n=\sum_{j=1}^{n-1}w_{n,j}x_{n-j}+r_n,\quad
w_{n,j}=a_j+\frac{b_j}{n}+c_{n,j}
\quad \lef n=1,2,\dots\ri\,,
\end{equation}
where $w_{n,j}\ge 0$, and $a,b,c$ are decaying at least exponentially
fast. The precise assumptions are formulated later. Our goal is to
prove that $x_n$ is polynomially decaying as $n\tv$ under suitable
conditions, and to determine the exponent.  

Similar recursions are widely examined, see e.g.\ Milne-Thompson
\cite{milne}, Cooper--Frieze \cite{cf}. In those cases either the
coefficients are special, or only the last $m$ terms appear on the
right-hand side for some fixed $m$. Now all previous terms are
present, and the weights depend both on $n$ and $j$.  

On the other hand, omitting $\lef b_n\ri$ and $\lef c_n\ri$ and
supposing that $\lef a_n\ri$ is a probability distribution, we get the
well-known renewal formula \cite[Chapter XIII]{feller1}. The
asymptotics were examined in a more subtle way in \cite{erdos}, for
instance, by weaker Tauberian type assumptions. In our case the
coefficients are of more general form; 
however, we have stronger conditions on them. An example will show
that these assumptions can not be totally omitted (see Remark
\ref{rekmj1}). 

The continuous counterpart is the following integral equation, which
is a Volterra equation of the second kind. 
\begin{equation}\label{inteq}
g(t)=\int_0^t w_{t,s}\, g(t-s)\,ds+r(t)
\end{equation}
for $t>0$ and $g\lef 0\ri=1$. The kernel $w_{t,s}$ is supposed to be
written in the following form.
\[
0\leq w_{t,s}=a(s)+\frac{b(s)}{t+d}+c_{t,s},
\]
where $a$ is a probability density function, and again, $a, b, c$ are
decreasing fast. We will show that $g(t)$ is between two
polynomially decaying functions under suitable conditions, and give
the exponent. In addition, assuming that 
$g$ is decreasing, we will prove that $g(t)$ is polynomially
decaying as $t\tv$. We use Laplace transforms and Tauberian theorems
in this part.   

Omitting $b$ and $c$ we get a classic renewal equation  \cite[Chapter
XI]{feller}.  

In Section 2 we formulate the main results for both cases. Sections 3
and 4 contain the proofs for the discrete and the continuous cases,
respectively. Section 5 contains the Laplace transform methods.

\section{Main results}

\subsection{The discrete recursion}

Consider the following recursion: 
\begin{equation}\label{reke1}
x_n=\sum_{j=1}^{n-1}w_{n,j}x_{n-j}+r_n,\quad
w_{n,j}=a_j+\frac{b_j}{n}+c_{n,j},
\quad \lef n=1,2,\dots\ri,
\end{equation}
where $w_{n,j}\ge 0$, and $a_n$, $b_n$, $c_{n,j}$, $r_n$ satisfy the 
following conditions.
\begin{itemize}
\item[(r1)] $a_n\ge 0$ for $n\geq 1$, and the greatest common
divisor of the set $\lk n:  a_n>0\rk$ is $1$;  
\item[(r2)] $r_n\ge 0$, and there exists such an $n$ that $r_n>0$; 
\item[(r3)] there exists $z>0$ such that
\begin{gather*}
1<\sum_{n=1}^\infty a_nz^n<\infty,\qquad
\sum_{n=1}^\infty|b_n|z^n<\infty,\\
\sum_{n=1}^\infty\sum_{j=1}^{n-1}|c_{n,j}|z^j<\infty,\qquad
\sum_{n=1}^\infty r_nz^n<\infty.
\end{gather*}
\end{itemize}
It is clear that $x_n\ge 0$ for $n\geq 1$.

Our theorem gives the asymptotics of $\lef x_n\ri$. It is 
polynomially decaying; the exponent is also given. 

\begin{theorem} \label{rekt1} Suppose that the sequence $(x_n)$
satisfies recursion \eqref{reke1}, conditions \textup{(r1)--(r3)} hold,
and $(x_n)$ has infinitely many positive terms. Then $x_n
n^{-\gamma}q^{n}\rightarrow C$ as $n\to\infty$, where $C$ is a
positive constant, $q$ is the positive solution of equation
$\sum_{n=1}^\infty a_nq^n=1$, and  
\[
\gamma=\frac{\sum\limits_{n=1}^\infty b_nq^n}
{\sum\limits_{n=1}^\infty na_nq^n}\,.
\] 
\end{theorem}

\begin{remark} The condition on  $w_{n,j}$ in recursion \eqref{reke1}
can be modified in the following way.
\[\label{reke2}
w_{n,j}=a_j+\frac{b_j}{n-j}+c_{n,j},\quad n=1,2,\dots\ .
\]
The difference may be added to the remainder term $c_{n,j}$, because
we have 
\begin{multline*}
\sum_{n=1}^\infty\sum_{j=1}^{n-1}\left|\frac{b_j}{n-j}-
\frac{b_j}{n}\right|y^j\\
=\sum_{j=1}^\infty|b_j|y^j\sum_{n=j+1}^\infty\Bigl(\frac{1}{n-j}-
\frac{1}{n}\Bigr)
=\sum_{j=1}^\infty|b_j|y^j\sum_{n=1}^j\frac 1n\,,
\end{multline*}
which is finite for $0<y<z$. Since the generating function of the
sequence  $(a_n)$ is left continuous at point $z$, there exists $y<z$
such that $\sum_{n=1}^\infty a_ny^n>1$.   
\end{remark}

\begin{remark} The condition that the sequence has infinitely many 
positive terms is necessary as the following example shows.  
Let  $r_1=1$,  $r_n=0$ if $n>1$, and
$w_{n,j}=a_n\Bigl(1-\dfrac{1}{n-j}\Bigr)$. Then we get that $x_1=1$,
$x_2=x_3=\dots=0$.
\end{remark}

\subsection{The integral equation}
Now we examine an integral equation, which is similar to recursion
\eqref{reke1}. Namely, let $g:\mathbb R\rightarrow \mathbb R$ be the
solution of the following integral equation; we will explain later why
the solution exists.  
\begin{equation}\label{rekie1}
g(t)=\int_0^t w_{t,s}\,g(t-s)\,ds+r(t) 
\end{equation}
for $t>0$ and $g\lef 0\ri=1$. Here 
\[
0\leq w_{t,s}=a(s)+\frac{b(s)}{t+d}+c_{t,s},
\]
and the following conditions hold. 
\begin{itemize}
\item[(i1)] $a\in L^1[0,\infty)$ is a probability density
function concentrated on the set of positive real numbers.  That is,
$a$ is nonnegative almost everywhere, and $\int_0^{\infty}a(s)\,ds=1$. 
\item[(i2)] $b\in L^1[0,\infty)$, and $d$ is a positive
constant. 
\item[(i3)] $r\in L^1[0,\infty)$ is a nonnegative, continuous
function.  
\item[(i4)] $c:\lsz 0,\infty\ri^2\rightarrow \R$ is (jointly)
measurable, $c_{t,s}$ is integrable on $\lsz 0,t\rsz$ with respect to
$s$ for all $t>0$, and $\lim_{t\to\infty}c_{t,s}=0$ for a.e. $s>0$. 
\item[(i5)] There exists $z>1$ such that 
\[
\int_0^{\infty} a(t) z^t\,dt<\infty,\quad \int_0^{\infty}
\left\vert b(t) \right\vert\,z^t\,dt<\infty,\quad \text{and}
\]
\item[(i6)]
$\displaystyle z^t\int_0^t|c_{t,s}|\,ds$ and $\displaystyle r(t) 
z^t$ are directly Riemann integrable with respect to $t$ on $\lsz
0,\infty\ri$. 
\end{itemize}

Recall that a nonnegative function $h$ is directly Riemann integrable
on $\lsz 0,\infty\ri$ (see p. 361 of \cite{feller}), if and only if it
is (Riemann) integrable on every finite interval, and for all $\tau>0$
we have
\[
\sum_{n=1}^{\infty} \sup_{n\tau\leq \theta\leq (n+1)\tau}
h(\theta)<\infty,
\]  
that is, the upper Riemann sum of $h$ with span $\tau$ is finite. As usual,
we say that a real function $h$ is directly Riemann integrable if both
its positive and negative parts are directly Riemann integrable. This
is equivalent to the direct Riemann integrability of $\vert h\vert$. 

Equation \eqref{rekie1} is a nonlinear Volterra-type integral equation
of the second kind. It is easy to check that all conditions of Theorem
3.2.\ of \cite{linz} hold for this equation in a finite interval
$0\leq t\leq T$. Thus, applying the theorem we get that the equation
has a unique and continuous solution for all positive $T$. Hence
$g(t)$ is defined on the set of nonnegative real numbers, and it is 
continuous. Since the proof of  Theorem 3.2. of \cite{linz} relies on
Picard approximation, and $g(0)$, $w$, $r$ are all nonnegative,
it is clear that $g(t)$ is nonnegative for all $t\geq 0$.  

Our main results are about the asymptotics of $g(t)$ as
$t\tv$. First we give the order of $g$ by proving lower and upper
bounds. Then assuming that $g$ is decreasing, we will find the
asymptotics of $g$ using Laplace transforms. 

\begin{theorem}\label{rekt2}
Let $g$ be the solution of equation \eqref{rekie1}. Suppose that $w$
is nonnegative, all conditions \textup{(i1)--(i6)} hold, and for
all $T>0$ there exists $t>T$ such that $g(t)>0$. Introduce 
\[
\gamma=\frac{\int_0^{\infty} b(s)\,ds}{\int_0^{\infty} sa(s)\,ds}.
\]  
Then $0<\liminf_{t\tv}g(t) t^{-\gamma}\leq \sup_{t}g(t)
t^{-\gamma}<\infty$ holds. 
\end{theorem}

\begin{theorem}\label{rekt3}
Let $g$ be the solution of equation \eqref{rekie1}. In addition to the
conditions of Theorem \ref{rekt2} suppose that $g$ is decreasing. Then
$g(t) t^{-\gamma}\rightarrow C$ holds for some $0<C<\infty$ as
$t\tv$.  
\end{theorem}

\section{The discrete case: Proof of Theorem 1.}

\subsection{Preliminaries} \label{rekpre}
We may assume that $\sum_{n=1}^{\infty} a_n=1$ and $q=1$.

Condition (r3) implies that  $q$ exists, and $q<z$. 
Define $\tilde
x_n=q^nx_n$, $\tilde a_j=q^ja_j$, $\tilde b_j=q^jb_j$, $\tilde
c_{n,j}=q^jc_{n,j}$, $\tilde r_n=q^nr_n$ for $n,j\geq 1$. We get that 
$\sum_{n=1}^\infty\tilde a_n=1$ and
\[
\tilde x_n=\sum_{j=1}^{n-1}\Bigl(\tilde a_j+\frac{\tilde b_j}{n-j}+
\tilde c_{n,j}\Bigr)\tilde x_{n-j}+\tilde r_n, \quad n=1,2,\ldots.
\]
Moreover, condition (r3) holds with $\tilde z=z/q$. Thus 
we may assume that $\sum_{n=1}^\infty a_n=1$, and $z>1$, indeed.

\begin{lemma} \label{rekl2} $x_n>0$ for every $n$ large enough.
\end{lemma}
\biz
If $a_k>0$ for some integer $k$, then for every sufficiently large $n$
we have $w_{n,k}>0$. To see this, note that
$\lim\limits_{n\to\infty}\frac{b_k}{n-k}=0$,  
and  $\lim\limits_{n\to\infty}
c_{n,k}=0$ holds for fixed $k$ due to condition (r3). Hence, if
$x_n>0$ and $n$ is large enough, then $x_{n+k}>0$. This  
implies that
$x_{n+\ell k}>0$ for $\ell=1,2,\dots$. Due to condition (r1), every
sufficiently large $n$ 
 is a linear combination of some values of $k$ for which
$a_k>0$. Therefore $x_n$ is positive for every $n$ large enough. \qed

We need some more notations.

Define $y_n=x_nn^{-\gamma}$ for $n\geq 1$. We have 
\begin{equation}\label{reke3}
y_n=\sum_{j=1}^{n-1}w_{n,j}\Bigl(1-\frac jn\Bigr)^{\gamma}y_{n-j}+
r_nn^{-\gamma}, \quad  n=1,2,\ldots\,.
\end{equation}
From the Taylor expansion of the function $f(x)=(1-x)^t$ for  $x\ge 0$ 
we get that
\begin{gather*}
\left|f(x)-1+tx\right|=\left|\tfrac{t(t-1)}2\,x^2(1-\theta
  x)^{t-2}\right|\le \tfrac{|t(t-1)|}{2}\,x^2\,e^{-\theta x(t-2)}\\
\le\tfrac{|t(t-1)|}{2}\,x^2\,e^{x|t-2|}\le
\tfrac{|t(t-1)|}{2}\,x^2\,e^{n\veps x},
\end{gather*}
where $0\le\theta\le 1$, and $\veps>0$ is so small that
 $e^{2\veps}<z$ holds with $z$ of condition (r3), while $n$ is so
 large that $n\veps>|t-2|$ holds.

Therefore
\begin{equation}\label{reke2.0}
\Bigl(1-\frac jn\Bigr)^{\gamma}=1-\frac{\gamma j}n+R_{n,j},
\end{equation}
 where 
\begin{equation}\label{reke2.1}
|R_{n,j}|\le\frac{|\gamma(\gamma-1)|}{2}\,\frac{j^2}{n^2}\,e^{j\veps} 
\end{equation}
holds uniformly in $j$ for all $n$ large enough, say $n\ge L$.
Assuming that the coefficients satisfy equation \eqref{reke1} we get
that   
\begin{equation}\label{reke4}
w_{n,j}\Bigl(1-\frac jn\Bigr)^{\gamma}=a_j
+\frac 1n\bigl(b_j-\gamma ja_j\bigr)
-\frac{1}{n^2}\,\gamma jb_j+w_{n,j}R_{n,j}+c_{n,j}
\Bigl(1-\frac{\gamma j}n\Bigr).
\end{equation}

\subsection{Boundedness of $\lef y_n\ri$}

Our next goal is to prove that the sequence $\lef y_n\ri$ is bounded
from above, and its limes inferior is positive. Before doing so
we prove another lemma. 

\begin{lemma} \label{rekl1}
For every positive integer $k$ we have
\begin{equation}
\sum_{n=k}^\infty\left|\sum_{j=1}^{n-k}w_{n,j}
\Bigl(1-\frac jn\Bigr)^\gamma-a_j\right|<\infty.\label{reke5}
\end{equation}
\end{lemma} 
\biz
Using equation \eqref{reke4} we obtain that
\begin{multline*}
\sum_{n=k}^\infty\left|\sum_{j=1}^{n-k}w_{n,j}
\Bigl(1-\frac jn\Bigr)^\gamma-a_j\right|=
\sum_{n=k}^\infty\biggl|\frac 1n\sum_{j=1}^{n-k}
\bigl(b_j-\gamma ja_j\bigr)-\frac{\gamma}{n^2}\sum_{j=1}^{n-k}jb_j
+\\
\qquad\qquad\qquad+\sum_{j=1}^{n-k}w_{n,j}R_{n,j}+
\sum_{j=1}^{n-k}c_{n,j}\Bigl(1-\frac{\gamma j}n\Bigr)\biggr|\\
\le \sum_{n=k}^\infty\Biggl(\frac 1n\sum_{j=n-k+1}^\infty
\bigl|\gamma ja_j-b_j\bigr|+\frac{|\gamma|}{n^2}\sum_{j=1}^{n-1}
j|b_j|+\\
\qquad\qquad\qquad+\sum_{j=1}^{n-1}|w_{n,j}R_{n,j}|+
\bigl(1+|\gamma|\bigr)\sum_{j=1}^{n-1}|c_{n,j}|\Biggr).
\end{multline*}

For the first term we used that 
$\sum_{j=1}^\infty\bigl(b_j-\gamma ja_j\bigr)=0$
holds by the definition of $\gamma$.

Let us divide the sum into four parts and examine them separately. By
condition (r3) for the first two we have that
\begin{gather*}
\sum_{n=k}^\infty\frac 1n\sum_{j=n-k+1}^\infty\bigl|\gamma
ja_j-b_j\bigr|
\le\sum_{n=k}^\infty\sum_{j=n-k+1}^\infty\bigl(|\gamma|
ja_j+|b_j|\bigr)\\=|\gamma|\sum_{j=1}^\infty j^2a_j+\sum_{j=1}^\infty j|b_j|
<\infty;\\
\sum_{n=1}^\infty\frac{|\gamma|}{n^2}\sum_{j=1}^{n-1}j|b_j|\le
|\gamma|\,\sum_{n=1}^\infty\frac{1}{n^2}\;\sum_{j=1}^\infty j|b_j|
<\infty.
\end{gather*}
 Moreover, using \eqref{reke2.1}, we obtain that
\begin{multline*}
\sum_{n=L}^\infty\sum_{j=1}^{n-1}|w_{n,j}R_{n,j}|\le
\frac{|\gamma(\gamma-1)|}2\,\sum_{n=L}^\infty\sum_{j=1}^{n-1}
\bigl(a_j+|b_j|+|c_{n,j}|\bigr)\frac{j^2}{n^2}e^{j\veps}\\
\leq\frac{|\gamma(\gamma-1)|}2\,\Bigg[\sum_{j=1}^\infty
\sum_{n=j+1}^{\infty}\bigl(a_j+|b_j|\bigr)\frac{j^2}{n^2}e^{j\veps}+  
\sum_{n=1}^\infty\sum_{j=1}^{n-1}|c_{n,j}|z^j\Bigg]<\infty
\end{multline*}
by condition (r3) and the choice of $\veps$.

Similarly,
\[
\sum_{n=1}^\infty\sum_{j=1}^{n-1}|c_{n,j}|\le
\sum_{n=1}^\infty\sum_{j=1}^{n-1}|c_{n,j}|z^j<\infty
\]
also holds. Thus we have proved \eqref{reke5}. \qed

\begin{lemma}\label{rekl3} 
$(y_n)$ is bounded from above.
\end{lemma}
\biz
Let $z_n=\max\bigl\{1,\,y_1,\,\dots,\,y_n\bigr\}$. Then $y_n\le z_n$,
$z_n$ is increasing, and
\begin{align*}
z_n&\le z_{n-1}\;\max\Biggl\{1,\
\sum_{j=1}^{n-1}w_{n,j}\Bigl(1-\frac jn\Bigr)^\gamma 
+r_nn^{-\gamma}\Biggr\}\\
&\le z_{n-1}\Biggl(1+\left|\sum_{j=1}^{n-1}
w_{n,j}\Bigl(1-\frac jn\Bigr)^\gamma-a_j\right|
+r_nn^{-\gamma}\Biggr)\\
&=z_{n-1}(1+s_n),
\end{align*}
where $s_n=\left|\sum_{j=1}^{n-1}
w_{n,j}\left(1-\frac jn\right)^\gamma-a_j\right|
+r_nn^{-\gamma}$.

Iterating this we obtain that $\sup_{n\ge 1}
z_n\le\prod_{n=1}^\infty (1+s_n)$. In  order to show that this
quantity is finite it is sufficient to prove that $\sum_{n=1}^\infty
s_n<\infty$ holds.  
The latter is implied by Lemma \ref{rekl1} with $k=1$, and by the
fact that
\[
\sum_{n=1}^\infty r_nn^{-\gamma}=
O\Biggl(\sum_{n=1}^\infty r_n z^n\Biggr)<\infty.
\]
Thus the sequence $\lef y_n\ri$ is bounded from above. \qed

\begin{lemma}\label{rekl4} 
$\liminf_{n\tv} y_n>0$.
\end{lemma}
\biz
This is similar to the upper bound, therefore we only outline the
proof, omitting the details.

Based on Lemma \ref{rekl2}, suppose that $x_n>0$ for all $n\ge N$. By
condition (r2) and the definition of $y_n$ we have that 
\[
y_n\ge \sum_{j=1}^{n-N}w_{n,j}\Bigl(1-\frac jn\Bigr)^\gamma y_{n-j}.
\]

We obtain that for  $z_n=\min\limits_{N\le j\le n}y_j$ the 
following inequality holds.
\begin{align*}
z_n&\ge z_{n-1}\,\min\Biggl\{1,\ \sum_{j=1}^{n-N}w_{n,j}
\Bigl(1-\frac jn\Bigr)^\gamma\Biggr\}\\  
&\ge z_{n-1}\Biggl(\ \sum_{j=1}^{n-N}a_j-\left|\sum_{j=1}^{n-N}
w_{n,j}\Bigl(1-\frac jn\Bigr)^\gamma-a_j\right|\,\Biggr)\\
&=z_{n-1}\Biggl(1-\sum_{j=n-N+1}^\infty a_j-\left|\sum_{j=1}^{n-N}
w_{n,j}\Bigl(1-\frac jn\Bigr)^\gamma-a_j\right|\,\Biggr)\\
&=z_{n-1}(1-s_n),
\end{align*}
where \[s_n=\sum_{j=n-N+1}^\infty a_j+\left|\sum_{j=1}^{n-N}
w_{n,j}\Bigl(1-\frac jn\Bigr)^\gamma-a_j\right|.\]

This implies that $\inf\limits_{n\ge N}y_n\ge\lim\limits_{n\to\infty}z_n
\ge z_N\prod_{n=N+1}^\infty(1-s_n)$. For the proof of the positivity
of the right-hand side we need to show that
$\sum_{n=N+1}^\infty s_n<\infty$. This is a consequence of Lemma
\ref{rekl1} with $k=N$, and the following estimation.
\[
\sum_{n=N+1}^\infty\ \sum_{j=n-N+1}^\infty a_j\le \sum_{j=1}^\infty
ja_j<\infty.
\]
We conclude that  $\inf\limits_{n\ge N}y_n>0$, and hence
$\liminf_{n\tv} y_n>0$. \qed 

\subsection{Final step}

The remaining part of the proof is similar to the proof of 
the discrete renewal theorem.

Equation \eqref{reke4} implies that
\begin{multline}\label{reke6}
y_n=\sum_{j=1}^{n-1}\biggl(a_j+\frac{b_j-\gamma ja_j}{n-j}-
\frac{(b_j-\gamma ja_j)j}{n(n-j)}-\frac{1}{n^2}\,\gamma jb_j+\\
+w_{n,j}R_{n,j}+c_{n,j}\Bigl(1-\frac{\gamma j}n\Bigr)\biggr)y_{n-j}+
r_nn^{-\gamma}.
\end{multline}

Fix a positive integer $N$. Then we obtain from \eqref{reke6} by
summation that 
\begin{equation}
\sum_{n=1}^Ny_n=\sum_{n=1}^N\sum_{j=1}^{n-1}\biggl(a_j
+\frac{b_j-\gamma ja_j}{n-j}\biggr)y_{n-j}+u_N\label{reke7}
\end{equation}
with an appropriately chosen sequence $\lef u_N\ri$. Here $u_N$ is
convergent as $N\to\infty$; let $u$ denote the limit. In order to show  
this, since $\lef y_n\ri$ is bounded, it is sufficient to prove that 
\begin{gather*}
\sum_{n=1}^\infty\sum_{j=1}^{n-1}\biggl(\frac{|b_j-\gamma ja_j|j}{n(n-j)}+
\frac{1}{n^2}\,|\gamma b_j|j
+|w_{n,j}R_{n,j}|+|c_{n,j}|\bigl(1+|\gamma|\bigr)\biggr)<\infty;\\
\sum_{n=1}^{\infty} r_nn^{-\gamma}<\infty.
\end{gather*}
We have almost done it before; the only thing left is to show the
convergence for the first term in the double sum. In this 
case 
\begin{multline*}
\sum_{n=1}^\infty\sum_{j=1}^{n-1}\frac{|b_j-\gamma ja_j|j}{n(n-j)}\le
\sum_{n=1}^\infty\sum_{j=1}^{n-1}\frac{j|b_j|+|\gamma| j^2a_j}
{n(n-j)}\\
=\sum_{j=1}^\infty\bigl(j^2|b_j|+|\gamma| j^3a_j\bigr)\sum_{n=j+1}^\infty
\frac{1}{n(n-j)}\\
\le
\sum_{j=1}^\infty\bigl(j^2|b_j|+|\gamma|
j^3a_j\bigr)\,\sum_{n=1}^\infty\frac 1{n^2} <\infty.
\end{multline*}

Introduce variable $m=n-j$ instead of $n$ in equation
\eqref{reke7}. Then 
\begin{align*}
\sum_{n=1}^Ny_n&=\sum_{j=1}^{N-1}\sum_{n=j+1}^N\biggl(a_j
+\frac{b_j-\gamma ja_j}{n-j}\biggr)y_{n-j}+u_N\\
&=\sum_{j=1}^{N-1}\sum_{m=1}^{N-j}\biggl(a_j
+\frac{b_j-\gamma ja_j}{m}\biggr)y_m+u_N\\
&=\sum_{m=1}^{N-1} y_m\sum_{j=1}^{N-m}\biggl(a_j
+\frac{b_j-\gamma ja_j}{m}\biggr)+u_N.
\end{align*}

Since
\[
\sum_{j=1}^\infty\biggl(a_j+\frac{b_j-\gamma ja_j}{m}\biggr)=1,
\]
we have
\[
\sum_{n=1}^Ny_n=\sum_{m=1}^N y_m\Biggl(1-\sum_{j=N-m+1}^\infty
\biggl(a_j +\frac{b_j-\gamma ja_j}{m}\biggr)\Biggr)+u_N.
\]
These imply that
\begin{multline*}
\sum_{m=0}^{N-1} y_{N-m}\sum_{j=m+1}^\infty a_j=
\sum_{m=1}^N y_m\sum_{j=N-m+1}^\infty a_j\\
=u_N-
\sum_{m=1}^N y_m\sum_{j=N-m+1}^\infty\frac{b_j-\gamma ja_j}{m}\,.
\end{multline*}
The second term on the right-hand side converges to $0$ as
$N\to\infty$, because 
\[
\sum_{m=1}^N y_m\sum_{j=N-m+1}^\infty\frac{|b_j-\gamma ja_j|}{m}
\le(\sup_n y_n)\sum_{m=1}^N\frac 1m\sum_{j=N-m+1}^\infty\bigl(
|b_j|+|\gamma|ja_j\bigr).
\]
The sum inside is estimated in the following way.
Let $\veps>0$ such that  $e^{\veps}<z$ holds. Then 
\[
\sum_{j=N-m+1}^\infty\bigl(|b_j|+|\gamma|ja_j\bigr)\le
K\,e^{-(N-m)\veps},
\]
where
\[
K=\sum_{j=1}^\infty\bigl(|b_j|+|\gamma|ja_j\bigr)e^{j\veps}<\infty.
\]
Now, using notation $M=\left[\sqrt{N}\right]$, we get that 
\begin{multline*}
\sum_{m=1}^N\frac Km\,e^{-(N-m)\veps}
\le \sum_{m=1}^{N-M}\frac Km\,e^{-(N-m)\veps}+
\sum_{m=N-M+1}^N \frac Km\,e^{-(N-m)\veps}\\
\le NKe^{-M\veps}+\frac{MK}{N-M}\,,
\end{multline*}
which tends to $0$ as $N\tv$. We conclude that
\begin{equation}
\sum_{m=0}^{N-1} y_{N-m}\sum_{j=m+1}^\infty a_j\to u\label{reke8}
\quad \lef N\tv\ri. 
\end{equation}

Modifying the proof of Lemma \ref{rekl1}, namely, using equation
\eqref{reke4}, condition (r3) and the fact that the sequence of
arithmetic means converges to zero if the original sequence is
nonnegative and converges to zero, it is easy to see that
\[
\sum_{j=1}^{n-1}\left|w_{n,j}
\Bigl(1-\frac jn\Bigr)^\gamma-a_j\right|\to 0 \quad \lef n\tv\ri.
\]
This and \eqref{reke3}, together with the boundedness of
$\lef y_n\ri$ imply that
\begin{equation}
y_n-\sum_{j=1}^{n-1}a_jy_{n-j}\to 0,\label{reke9}
\end{equation}
as $n\to\infty$.
 
From now on the argument is the usual one.

Let $\lef n_k\ri$ be a subsequence of the natural numbers that
satisfies 
\[
\lim_{k\to\infty}y_{n_k}=\limsup_{n\to\infty}y_n=:\overline y.
\]
From \eqref{reke9} we get that for all $\ell\le M$ the following
estimation holds.
\begin{align*}
\overline y=&\lim_{k\to\infty}y_{n_k}\le
a_\ell\liminf_{k\to\infty}y_{(n_k-\ell)}+ 
\limsup_{k\to\infty}\sum_{j<n_k,\,j\ne\ell}a_jy_{(n_k-j)}\\
&\le a_\ell\liminf_{k\to\infty}y_{(n_k-\ell)}+
\sum_{j\ne\ell,\,j\le M} a_j\limsup_{k\to\infty}y_{(n_k-j)}
+(\sup_n y_n)\sum_{j=M+1}^\infty a_j\\
&\le a_\ell\liminf_{k\to\infty}y_{(n_k-\ell)}+
(1-a_\ell)\overline y+(\sup_n y_n)\sum_{j=M+1}^\infty a_j.
\end{align*}

Since $M$ may be arbitrarily large, this immediately implies that 
 $\lim\limits_{k\to\infty} y_{(n_k-\ell)}=\overline y$,
for all $a_\ell>0$. By iteration we obtain that 
\[
\lim_{k\to\infty}y_{(n_k-\ell_1-\ldots-\ell_i)}=\overline y,
\]
for all positive $a_{\ell_1},\,\dots,\,a_{\ell_i}$. By
condition (r1) for all sufficiently large $m$ we have
$\lim\limits_{k\to\infty} y_{n_k-m}=\overline y$. Modifying the 
subsequence we may assume that this holds for all
$m=0,1,\dots$. Hence choosing $N=n_k$ in \eqref{reke8} we can see that 
\[
\overline y\,\sum_{m=0}^\infty\sum_{j=m+1}^\infty a_j= u.
\]

For $\underline y=\liminf\limits_{n\to\infty}y_n$  the same argument
shows that 
\[
\underline y\,\sum_{m=0}^\infty\sum_{j=m+1}^\infty a_j=u.
\]
Hence $\overline y=\underline y$, that is, the limit 
$\lim_{n\to\infty}y_n=C$ exists. We have already proved that this is 
finite and positive. This implies Theorem \ref{rekt1}.
\qed

\begin{remark}
It is well known that if $y_n=\sum_{j=1}^{n-1}a_jy_{n-j}$ holds for
all $n$, instead of in \eqref{reke9}, then $y_n$
is convergent. On the other hand, \eqref{reke9} is not yet sufficient
for the convergence of $\lef y_n\ri$.  
For example, let $y_n=2+\sin\bigl(\log(1+n)\bigr)$. Then
\[
\left|y_n-y_{n-j}\right|\le\log(1+n)-\log(1+n-j)\le
\frac{j}{1+n-j}\,,
\]
hence
\begin{align*}
\biggl|y_n-\sum_{j=1}^{n-1}a_j y_{n-j}\biggr|&\le
y_n \sum_{j=n}^\infty a_j+
\sum_{j=1}^{n-1}\left|y_n-y_{n-j}\right|a_j\\
&\le 3\sum_{j=n}^\infty a_j+3\sum_{j=1}^{n-1}\frac{ja_j}{1+n-j}\\
&\le 3\sum_{j=n}^\infty a_j+3\sum_{j=1}^{M-1}\frac{ja_j}{n-M}+
3\sum_{j=M}^{n-1}ja_j\\
&\le 3\sum_{j=n}^\infty a_j+\frac{3}{n-M}\sum_{j=1}^\infty ja_j+
3\sum_{j=M}^{\infty}ja_j.
\end{align*}
This converges to zero with $M=n/2$, but $y_n$ does not converge.
\end{remark}

\begin{remark} \label{rekmj1}
 If $w_{k,i}=a_i$, then $x_k\to C$ by the arithmetic version of the
 renewal theorem. The following example shows that a remainder,
 though converging to $0$, may change this. Let $(a_i)$ be arbitrary,
 $x_k=2+\sin\bigl(\log(k+1)\bigr)$, and 
\begin{equation*}
w_{k,i}=a_i+\biggl(x_k-\sum_{j=1}^{k-1}x_{k-j}a_j\biggr)
\biggl(\sum_{j=1}^{k-1}x_j\biggr)^{\!-1}.
\end{equation*}
Then $w_{k,i}=a_i+o(1)$, and
$x_k=\sum\limits_{i=1}^{k-1}x_{k-i}w_{k,i}+2\delta_{k,1}$.
\end{remark}

\section{The continuous case: Proof of Theorem 2.}

\subsection{Preliminaries}

\begin{lemma} \label{rekli1}
Under the conditions of Theorem \ref{rekt2} we have
$g(t)>0$ for every $t$ large enough.
\end{lemma} 
\biz
Choose $0<s_0$ and $\delta>0$ such that the set $S'=\{
s\in (0, s_0): a(s)>3\delta\}$ has positive
Lebesgue measure.

Let $\ell(s)=\sup\{t\ge s: |c_{t,s}|\ge\delta\}$. Then
$\ell(s)<\infty$ for a.e. $s>0$ by condition (i4). Though $\ell$ is
not necessarily Borel measurable, yet it is Lebesgue measurable by the
measurable projection theorem, for the superlevel set $\{\ell>K\}$ is 
just the projection of the two dimensional measurable set
$\{(s,t):0<s\le t,\ K<t,\ |c_{t,s}|\ge\delta\}$ onto the first
coordinate. Hence $U_t=\{\ell<t\}$, $t>0$ is an increasing family of
Lebesgue measurable sets, and the Lebesgue measure of
$(0,s_0)\setminus U_t$ tends to $0$ as $t\to\infty$. The same holds
for the sequence $V_t=\{s\in(0,t]:|b(s)|\le\delta(t+d)\}$. Thus we
can find a threshold $T\ge s_0$ such that the Lebesgue measure of
$S=S'\cap U_T\cap V_T$ is positive. Obviously,
$w_{t,s}\ge\delta$ for all $t\ge T$ and $s\in S$. By the Lebesgue
density theorem we may assume that $S$ only consists of points with
density $1$. 

By the continuity of $g$ there exists a whole open interval $I$ above
$T$ where $g(t)$ is separated from zero. Let $\varepsilon$ denote
the length of $I$, and $\eta>0$ the infimum of $g$ over $I$. Then
for $t\in S+I$ we have
\begin{multline*}
g(t)\ge\int_0^t w_{t,s}\,g(t-s)\,ds\ge
\int_{S\cap(t-I)}w_{t,s}\,g(t-s)\,ds\\
\ge \delta\eta\,\lambda\bigl(S\cap(t-I)\bigr)>0,
\end{multline*} 
where $\lambda$ stands for the Lebesgue measure.

Since the set sum $S+I$ is an open set, we can iterate this procedure
to obtain that $g$ is positive everywhere on the set
\[
I\cup(S+I)\cup(S+S+I)\cup(S+S+S+I)\cup\dots\,.
\]

The proof can be completed by showing that this set contains every
sufficiently large real number. In other words, if $t$ is large
enough, then it can be written in the form $t=s_1+\dots +s_n+r$,
where $s_1,\dots,s_n\in S$, $n\in\mathbb N$, and $0<r<\varepsilon$.
In fact, this is true for arbitrary $S\subset\mathbb R^+$ that has
two incommensurable elements $\alpha$ and $\beta$ (hence for every set
of positive Lebesgue measure). Indeed, by the equidistribution theorem
there exists positive integers $k$ and $m$ such that
\[
k< m\,\frac{\alpha}{\beta}<k+\frac{\varepsilon}{\beta}\,,
\]
that is, $k\beta<m\alpha$, and their distance is less than
$\varepsilon$. Consequently, in the finite sequence
\[
nk\beta<(n-1)k\beta+m\alpha<(n-2)k\beta+2m\alpha<\dots<nm\alpha
\]
the distance between neighbouring terms is less than $\varepsilon$. If
$n$ is large enough, then $(n+1)k\beta<nm\alpha$, i.e., the largest term
of the sequence above is bigger than the smallest term of the next
sequence. Thus every $t\ge nk\beta$ is sufficiently close to a
positive linear combination of $\alpha$ and $\beta$.\qed

Let us introduce the notation 
\[
H(t)=g(t)\lef t+d\ri^{-\gamma} \qquad \lef t\geq 0\ri.
\]
From \eqref{rekie1} we obtain the following integral equation for $H$. 
\begin{equation}\label{rekie2}
 H(t)=\int_0^t w_{t,s}\lef \frac{t-s+d}{t+d}\ri^{\gamma}H(t-s)\,
 ds+r(t)\lef t+d\ri^{-\gamma} 
\end{equation}
for $t>0$, and $H\lef 0\ri=d^{-\gamma}$.

Let us choose $\varepsilon$ in a similar way as we did in the discrete
case. Namely, we have $e^{2\veps}<z$ with $z$ of condition (i5). In
what follows equations \eqref{rekie2.0}, \eqref{rekie2.1}, and
\eqref{rekie4} correspond to \eqref{reke2.0}, \eqref{reke2.1}, and
\eqref{reke4}, resp. Firstly,
\begin{equation}\label{rekie2.0}
\lef \frac{t-s+d}{t+d}\ri^{\gamma}=\lef
1-\frac{s}{t+d}\ri^{\gamma}=1-\frac{\gamma s}{t+d}+R_{t,s},
\end{equation}
where 
\begin{equation}\label{rekie2.1}
\left\vert R_{t,s}\right\vert\leq \frac{\left\vert \gamma
(\gamma-1)\right\vert}{2}\,\frac{s^2}{t^2}\,e^{s\veps},
\end{equation}
if $t$ is large enough, say $t\ge L$. Finally, from the decomposition
of $w_{t,s}$ we get that   
\begin{multline}\label{rekie4}
w_{t,s}\Bigl(1-\frac s{t+d}\Bigr)^{\gamma}\\
=a(s)+\frac {b(s)-\gamma sa(s)}{t+d}
-\frac{\gamma sb(s)}{(t+d)^2}+w_{t,s}R_{t,s}+c_{t,s}
\Bigl(1-\frac{\gamma s}{t+d}\Bigr).
\end{multline}

\subsection{Boundedness of $H$} The method of proof is
discretization; in this way all we need to do is similar to what we
did in the discrete case.

Before proving boundedness we need another lemma.

\begin{lemma}\label{rekli2}
For every fixed $T\geq0$ the function
\[
A(t)=\left\vert\int_0^{t-T}\lef w_{t,s}\Bigl(1-\frac{s}{t+d}\Bigr)
^{\gamma}-a(s)\ri ds\right\vert
\]
is directly Riemann integrable on $\left[T,\infty\ri$.
\end{lemma}
\biz
Fix $\tau>0$. Then 
\begin{equation}\label{rekie4.0}
\sum_{n=1}^{\infty}\sup_{n\tau\leq\theta\leq(n+1)\tau}
\int_0^{\theta}\left\vert c_{\theta,s}\right\vert z^s ds<\infty, 
\end{equation}
according to condition (i6). Now we prove that
\begin{equation}\label{rekie4.1}
\sum_{n=\bigl\lceil\tfrac{T}{\tau}\bigr\rceil}^{\infty}\ 
\sup_{n\tau\leq\theta\leq (n+1)\tau} A(\theta)=
\sum_{n=\bigl\lceil\tfrac{T}{\tau}\bigr\rceil}^{\infty}\ 
\sup_{0\leq\theta\leq\tau}A(n\tau+\theta)<\infty.
\end{equation} 
From the definition of $\gamma$ it follows that $\int_0^{\infty}
\bigl(b(s)-\gamma s a(s)\bigr) ds=0$. Using this and equation
\eqref{rekie4} we obtain that   
\begin{gather*}
A(n\tau+\theta)=\Bigg\vert\int_0^{n\tau+\theta-T}\lef
w_{n\tau+\theta, s}\Bigl(1-\frac{s}{n\tau+\theta+d}\Bigr)^{\gamma}-a(s)
\ri ds\Bigg\vert \\
= \Bigg\vert\int_0^{n\tau+\theta-T}\frac{ b(s)-\gamma s
 a(s)}{n\tau+\theta+d}\,ds-\frac{\gamma}{(n\tau+\theta+d)^2}
\int_0^{n\tau+\theta-T} sb(s)\,ds\\ 
+\int_0^{n\tau+\theta-T}w_{n\tau+\theta,s}
 R_{n\tau+\theta,s}\,ds+\int_0^{n\tau+\theta-T}c_{n\tau+\theta,s}
\Bigl(1-\frac{\gamma s}{n\tau+\theta+d}\Bigr)ds\Bigg\vert\\ 
\le\int_{n\tau+\theta-T}^{\infty}\frac{|b(s)-\gamma s a(s)|}
{n\tau+\theta+d}\, ds+\frac{|\gamma|}{(n\tau+\theta+d)^2}
\int_0^{n\tau+\theta} s\left\vert b(s)\right\vert ds\\ 
+\int_0^{n\tau+\theta}\left\vert w_{n\tau+\theta,s}
R_{n\tau+\theta,s}\right\vert ds+\bigl(1+|\gamma|\bigr)
\int_0^{n\tau+\theta}\left\vert
 c_{n\tau+\theta,s}\right\vert ds 
\end{gather*}
holds for all $\theta>0$. 

We treat the four integrals in the right-hand side separately again.
For the first term we get that 
\begin{align*}
\sum_{n=\big\lceil\tfrac{T}{\tau}\big\rceil}^{\infty}\
\sup_{0\leq\theta\leq\tau}\frac{1}{n\tau+\theta+d}
&\int_{n\tau+\theta-T}^{\infty}|b(s)-\gamma s a(s)|\,ds\\
&\leq \sum_{n=\big\lceil\tfrac{T}{\tau}\big\rceil}^{\infty}\
\frac{1}{n\tau+d}\int_{n\tau-T}^{\infty}|b(s)-\gamma s a(s)|\,ds\\ 
&\leq \sum_{n=1}^{\infty}\frac{n}{n\tau+d+T}
\int_{n\tau}^{(n+1)\tau}|b(s)-\gamma s a(s)|\,ds\\ 
&\leq \frac{1}{\tau} \int_{0}^{\infty}\bigl(|\gamma|s a(s)+|b(s)|
\bigr)ds,
\end{align*}
which is finite by condition (i5).

The sum of the second terms is also finite by condition (i5). 
\begin{multline*}
\sum_{n=1}^{\infty}\ \sup_{0\leq\theta\leq\tau}\frac{|\gamma|}
{(n\tau+\theta+d)^2}\int_0^{n\tau+\theta} s\left|b(s)\right| ds\\
\leq\frac{|\gamma|}{\tau^2}\,\sum_{n=1}^{\infty}\,\frac{1}{n^2}
\int_{0}^{\infty}s\left\vert b(s)\right\vert ds<\infty.
\end{multline*}
 
By inequality \eqref{rekie2.1}, condition $(i5)$, and the choice of
$\varepsilon$ in the preliminaries, we obtain the following estimation
for the sum of the third terms. 
\begin{multline*}
\sum_{n=\big\lceil\tfrac{L}{\tau}\big\rceil}^{\infty}
\sup_{0\leq \theta\leq\tau}\int_0^{n\tau+\theta}\left|w_{n\tau+\theta,
s}\,R_{n\tau+\theta,s}\right|ds\\ 
\leq\sum_{n=\big\lceil\tfrac{L}{\tau}\big\rceil}^{\infty}
\sup_{0\leq\theta\leq\tau}\int_0^{n\tau+\theta}\lef
a(s)+\frac{|b(s)|}d+|c_{n\tau+\theta,s}|\ri \frac{\left|\gamma
(\gamma-1)\right|s^2}{2(n\tau+\theta)^2}\,e^{s\varepsilon}\\ 
\leq \frac{\left|\gamma(\gamma-1)\right|}{2}\,\sum_{n=1}^{\infty} 
\int_0^{(n+1)\tau} \lef a(s)+\frac{|b(s)|}{d}\ri \frac{s^2}{\tau^2n^2}
\,e^{s\varepsilon}\,ds\\ 
+\frac{\left|\gamma(\gamma-1)\right|}{2}\,\sum_{n=1}^{\infty}\  
\sup_{0\leq\theta\leq\tau}\int_{0}^{n\tau+\theta} 
\left|c_{n\tau+\theta,s}\right|z^s\,ds. 
\end{multline*}
In the right-hand side the last sum is finite by \eqref{rekie4.0}. In
the first sum the integrand is nonnegative, hence by Fubini's theorem
we obtain that  
\begin{multline*}
\sum_{n=1}^{\infty} \int_0^{(n+1)\tau}\lef a(s)+\frac{|b(s)|}{d}\ri 
\frac{s^2}{\tau^2n^2}\,e^{s\varepsilon}\,ds \\
\leq \frac{1}{\tau^2}\int_0^{\infty} \lef\sum_{n=1}^{\infty}
\frac{1}{n^2} \ri\lef a(s)+\frac{\left\vert b(s)\right\vert}{d}\ri s^2
e^{s\varepsilon}\,ds<\infty, 
\end{multline*}
by the choice of $\varepsilon$. Thus the sum of the third terms is
finite, too. 

Finally, for the sum of the fourth terms we clearly have 
\[
\sum_{n=1}^{\infty} \sup_{0\leq \theta\leq\tau}
\int_0^{n\tau+\theta}\left|c_{n\tau+\theta,s}\right|ds\leq  
\sum_{n=1}^{\infty} \sup_{0\leq \theta\leq\tau}
z^{n\tau+\theta}\int_0^{n\tau+\theta} \left|
c_{n\tau+\theta,s}\right| ds<\infty.
\]

Putting these together we obtain that all four parts of $A$ give
finite sums, hence \eqref{rekie4.1} holds. The nonnegativity and
integrability of $A$ is clear, for it is a continuous function of
$t$. Thus the proof of the lemma is completed. \qed 

\begin{lemma}\label{rekli3}
$H(t)$ is bounded from above.
\end{lemma}

We define $Z(t)=\max\left\lbrace 1,\, H(s): 0\leq s\leq t\right\rbrace$
for $t\geq 0$. This is finite, because $H$, as well as $g$, is
continuous. 

First we give an upper bound for $\sup_{0<\theta\le\tau}H(t+\theta)$,
where $t$ and $\tau$ are fixed positive numbers. Introduce  
\[
w\cs_{t,s}=w_{t,s}\lef 1-\frac{s}{t+d}\ri^{\!\gamma}\qquad 0<t,\ 0\leq
s\leq t.
\] 
Using the nonnegativity of $w$ and $r$, equation \eqref{rekie2}, and
the definition of $Z$, we get that 
\begin{multline}\label{rekie4.2}
H(t+\theta)=\int_0^{\theta}w\cs_{t+\theta,s}H(t+\theta-s)\,ds\\
+\int_{\theta}^{t+\theta}w\cs_{t+\theta,s}H(t+\theta-s)\,ds
+r(t+\theta)(t+\theta+d)^{-\gamma}
\end{multline}
\[
\leq \int_0^{\theta}w\cs_{t+\theta,s}\,ds\,Z(t+\theta)
+\left[ \int_{\theta}^{t+\theta}w\cs_{t+\theta,s}\,ds+r(t+\theta)
(t+\theta+d)^{-\gamma}\right]Z(t)
\]
\begin{multline*}
=\int_0^{\theta}w\cs_{t+\theta,s}\,ds\,\lsz Z(t+\theta)-Z(t)\rsz
\\
+\left[\int_{0}^{t+\theta}w\cs_{t+\theta,s}\,ds+r(t+\theta)
(t+\theta+d)^{-\gamma}\right]Z(t). 
\end{multline*}

Next we want to prove that there exists $\tau_0>0$, and for every $\tau$,
$0<\tau\le\tau_0$, a positive integer $N(\tau)$ such that
\begin{equation}\label{rekie5}
\sup_{0\leq \theta\leq \tau} \int_0^{\theta}w_{n\tau+\theta,s}^*
\,ds\leq \frac{1}{2}\,,
\end{equation}
provided  $n>N(\tau)$.

To show this we will give an upper bound on
\begin{equation}\label{rekie5.1}
w_{n\tau+\theta,s}^*=\lef 1-\frac{s}{n\tau+\theta+d}\ri^{\!\gamma}\lef
a(s)+\frac{b(s)}{n\tau+\theta+d}+c_{n\tau+\theta,s}\ri.
\end{equation}

We clearly have
\[
\lef 1-\frac{s}{n\tau+\theta+d}\ri^{\!\gamma}\leq 
\lef 1+\frac{\theta}{d}\ri^{\!|\gamma|}\leq 
\exp\Big(\frac{\theta|\gamma|}{d}\Big),
\]
and
\[
a(s)+\frac{b(s)}{n\tau+\theta+d}\leq a(s)+\frac{|b(s)|}{d}\,.
\]
Hence, for $0\le\theta\le\tau$ we can write
\begin{multline*}
\int_0^{\theta}w_{n\tau+\theta,s}^*\,ds\le 
\exp\Big(\frac{\theta|\gamma|}{d}\Big)\left[\,
\int_0^\theta\Bigl(a(s)+\frac{|b(s)|}{d}\Bigr)ds+
\int_0^{\theta}\left|c_{n\tau+\theta,s}\right|ds\,\right]\\
\le\exp\Big(\frac{\tau|\gamma|}{d}\Big)
\int_0^\tau\Bigl(a(s)+\frac{|b(s)|}{d}\Bigr)ds+
z^{n\tau+\theta}\int_0^{\theta}\left|c_{n\tau+\theta,s}\right|ds,
\end{multline*}
if $n$ is large enough, namely, $n\ge |\gamma|/2d\varepsilon$ will do.
The first term in the right-hand side can be arbitrarily small if $\tau$
is fixed small enough. As to the second term, it can be estimated in the
following way. 
\[
\sup_{0\le\theta\le\tau}z^{n\tau+\theta}\int_0^{\theta}
\left|c_{n\tau+\theta,s}\right|ds\le
\sup_{n\tau\le\theta\le(n+1)\tau}z^{\theta}\int_0^{\theta}
\left|c_{\theta,s}\right|ds,
\]
which tends to $0$ as $n\to\infty$ by condition (i6). Thus
\eqref{rekie5} is satisfied if $n$ is greater than a certain threshold
$N(\tau)$. 

For any $0<\tau\leq \tau_0$ and $t=n\tau$, $n>N\lef \tau\ri$ inequality
\eqref{rekie4.2} implies that  
\begin{multline*}
\sup_{0< \theta \leq\tau} H(t+\theta)\leq 
\frac{1}{2}\lsz Z(t+\tau)-Z(t)\rsz+\\
+\sup_{0<\theta\leq\tau} \lsz\int_0^t w_{t+\theta,s}^*
ds+r(t+\theta+d)(t+\theta+d)^{-\gamma}\rsz Z(t).
\end{multline*}
Here we use that $Z$ is nonnegative and increasing by definition. 

We clearly have $Z\lef t+\tau\ri=\max\lk Z(t),\, \sup_{0<\theta\leq\tau}
H(t+\theta)\rk$. Therefore we obtain that  
\begin{multline*}
Z(t+\tau)\leq \max\Biggl\{Z(t),\quad \frac12\lsz Z(t+\tau)-Z(t)\rsz+\\
+\sup_{0<\theta\leq\tau}\left[\int_0^t w_{t+\theta,s}^*\,ds+
r(t+\theta)(t+\theta+d)^{-\gamma}\right]Z(t)\Biggr\},
\end{multline*}
from which it follows that
\begin{multline*}Z\lef t+\tau\ri-Z(t)\leq 
\Biggl(\frac12 \lsz Z\lef t+\tau\ri-Z(t)\rsz+\\
+\sup_{0<\theta\leq\tau}\lsz\int_0^tw_{t+\theta,s}^*\,ds-1+r(t+\theta)  
(t+\theta+d)^{-\gamma}\rsz Z(t)\Biggr)^{\!+};
\end{multline*}
where $x^+$ denotes $\max\lef x, 0\ri$, as usual. Hence
\begin{multline*}
Z(t+\tau)-Z(t)\\ 
\leq 2\Biggl(\sup_{0<\theta\leq\tau}\lsz\int_0^t w_{t+\theta,s}^*\,ds
-1+r(t+\theta)(t+\theta+d)^{-\gamma}\rsz\Biggr)^{\!+}Z(t).
\end{multline*}

We continue with deriving an upper bound for the right-hand side.
Since  $a$ is a probability density function, we have
\begin{multline*}
\sup_{0<\theta\leq\tau}\lsz\int_0^{t+\theta}w_{t+\theta,s}^*\,ds-1\rsz 
\le\sup_{0<\theta\leq\tau}\lsz\int_0^{t+\theta}w_{t+\theta,s}^*\,ds-
\int_0^{t+\theta}a(s)\,ds\rsz\\
\le\sup_{0<\theta\leq\tau}\left|\int_0^{t+\theta}
\big(w_{t+\theta,s}^*-a(s)\big)\,ds\right|.
\end{multline*}
Therefore  
\begin{multline*}
Z(t+\tau)-Z(t)\le Z(t)\Biggl(\sup_{0<\theta\leq\tau}
\left\vert \int_0^{t+\theta}\bigl(w_{t+\theta,s}^*-a(s)\bigr)ds
\right\vert+\\
+\sup_{0<\theta\leq\tau} r(t+\theta)(t+\theta+d)^{-\gamma}\Biggr)
\end{multline*}
for all $0<\tau\leq \tau_0$, $t=n\tau$, $n\geq N(\tau)$. 

Similarly to Lemma \ref{rekl3} of the discrete case, for the
boundedness of $Z(n\tau)$ from above it suffices to prove that   
\begin{multline}\label{rekie6}
\sum_{n=1}^{\infty}\,\sup_{0<\theta\leq\tau} \left|
\int_0^{n\tau+\theta}\bigl(w_{n\tau+\theta,s}^*-a(s)\bigr)ds\right|
+\\
+\sum_{n=1}^{\infty}\,\sup_{0\leq \theta\leq\tau}r(n\tau+\theta)
(n\tau+\theta+d)^{-\gamma} <\infty.
\end{multline}

Lemma \ref{rekli2} with $T=0$ implies that the first sum is finite.  
Since by condition (i6)  $r(t)z^t$ is directly Riemann integrable,
it follows that the second sum is also finite. Thus we conclude that
the sequence 
\[
Z(n\tau)=\max\lk 1,\,H(s): 0\leq s\leq n\tau\rk
\]
is bounded from above if $\tau$ is small enough. Hence the function
$H$ is also bounded from above.  \qed

\begin{lemma}\label{rekli4}
$\liminf_{t\tv}H(t)>0$.
\end{lemma}
\biz
Like in the discrete case, we omit the details that are
straightforward modifications of the previous lemma, and only give a
sketch of the proof. 

Based on Lemma \ref{rekli1}, we can suppose that $H(t)>0$ for all
$t\geq T$. This time define   
\[
Z(t)=\min\lk H(s): T\leq s\leq t\rk, 
\]
for $t\geq T$. 

Let us derive a lower bound for $H(t+\theta)$, where $t>T$ and
$\theta>0$. 
\begin{align*}
H(t&+\theta)=\int_{0}^{t+\theta} w^*_{t+\theta,s}\,H(t+\theta-s)\,ds+
r(t+\theta)(t+\theta+d)^{-\gamma}\\ 
&\ge\int_{0}^{\theta} w^*_{t+\theta,s}\,H(t+\theta-s)\,ds+
\int_{\theta}^{t-T+\theta} w^*_{t+\theta,s}\,H(t+\theta-s)\,ds\\
&\ge Z(t+\theta)\int_{0}^{\theta} w^*_{t+\theta,s}\,ds+
Z(t)\int_{\theta}^{t-T+\theta} w^*_{t+\theta,s}\,ds\\
&=\lsz Z(t+\theta)-Z(t)\rsz\int_{0}^{\theta} w^*_{t+\theta,s}\,ds+ 
Z(t)\int_{0}^{t-T+\theta} w^*_{t+\theta,s}\,ds.
\end{align*}

Now $Z$ is decreasing. Applying \eqref{rekie5} we obtain that  
\begin{equation*}
H(t+\theta)\geq \frac12 \lsz Z(t+\theta)-Z(t)\rsz
+Z(t) \int_{0}^{t-T+\theta} w^*_{t+\theta,s}\,ds
\end{equation*}
for $0<\theta\le\tau$. 
Taking infimum, subtracting $Z(t)$ and using that $a$ is a probability
density function we get that  
\begin{multline*}
Z(t+\tau)-Z(t)\ge\min\Bigg\{0,\quad\frac12\lsz Z(t+\tau)-Z(t)\rsz+\\
+Z(t)\lsz\inf_{0<\theta\leq \tau}\int_{0}^{t-T+\theta} 
w^*_{t+\theta,s}\,ds-1\rsz\Bigg\},
\end{multline*}
from which it follows that
\begin{multline*}
Z(t+\tau)-Z(t)\ge 2 \min\Bigg\{0,\ Z(t)\lsz\inf_{0<\theta\le\tau}
\int_{0}^{t-T+\theta} w^*_{t+\theta,s}\,ds-1\rsz \Bigg\}\\
\geq -2 \Biggl[\sup_{0<\theta\le\tau} \left|\int_{0}^{t-T+\theta} 
\bigl(w^*_{t+\theta,s}-a(s)\bigr)ds\right|+\int_{t-T}^{\infty} 
a(s)\,ds\Biggr] Z(t).
\end{multline*}

Similarly to Lemma \ref{rekl4}, in order to prove that $\lim_{n\tv} 
Z(T+n\tau)>0$ it suffices to show that  
\[
\sum_{n=1}^{\infty}\ \int_{(n-1)\tau}^{\infty} a(s)\,ds
+\sum_{n=1}^{\infty}\ 
\sup_{0<\theta\le\tau}\left|\int_{0}^{n\tau+\theta}
\bigl(w^*_{T+n\tau+\theta,s}-a(s)\bigr)ds\right|<\infty.
\]
For the first term we have
\[
\sum_{n=1}^{\infty}\ 
\int_{(n-1)\tau}^{\infty}a(s)\,ds\leq \frac{1}{\tau}\int_0^{\infty}
s a(s)\,ds<\infty
\]
by condition (i5). The finiteness of the second term follows directly
from Lemma \ref{rekli2}. 

Thus we proved that $\lim_{n\tv} Z(T+n\tau)>0$, which immediately
implies that $\liminf_{t\tv} H(t)>0$, as needed.\qed

\section{The monotonic case: Proof of Theorem 3}

First note that $\gamma\leq 0$ follows from the assumption that $g$ is
decreasing. 

We consider the following integral equation. 
\begin{multline}\label{reklap1}
g(x)=\int_0^xg(x-u)a(u)\,du+\int_0^xg(x-u)\,\frac{b(u)}{x+d}\,du + \\
+\int_0^x g(x-u)c_{x,u}\,du+r(x).  
\end{multline}

In the sequel we define the Laplace transform of an integrable
function $f$ as  
\[
F(s)=\lim_{y\tv} \int_0^y e^{-sx}f(x)\,dx
\]
for $s\in \mathbb C$, provided the limit exists and it is finite. 

Denote the Laplace transforms of functions $g$, $a$, $b$, $r$ by $G$,
$A$, $B$, $R$, respectively. These Laplace transforms are well defined
and holomorphic on the half-plane $\mathbb H=\lk s\in\mathbb C:
\text{Re } s>\rk$; this follows from the conditions on $a$, $b$, $r$
and Theorem \ref{rekt2}. Moreover, $A$, $B$, and $R$ are also holomorphic
in a neighbourhood of the origin. Let $f$ be either of the functions
above, then we have   
\[
F'(s)=-\int_0^{\infty} e^{-sx}xf(x)\,dx,\quad s\in\mathbb H.
\]

Multiplying both sides of equation \eqref{reklap1} by $(x+d)e^{-sx}$,
then integrating, and using the well known properties of the Laplace
transform we obtain that
\begin{align*}
-G'(s)&+d\cdot G(s)=-\lsz G(s) A(s)\rsz'+\\
&+d\cdot G(s) A(s)+G(s) B(s)+C(s)-R'(s)+d\cdot R(s)
\end{align*}
for $s\in \mathbb H$, where
\[
C(s)=\int_0^\infty e^{-sx}(x+d)\int_0^x g(x-u)c_{x,u}\,du\,dx. 
\]
This is finite and holomorphic in a neighbourhood of $0$ by condition
(i6) and Theorem \ref{rekt2}.

After rearranging we have
\begin{align*}
G'(s)&=G(s)\lsz d-\frac{B(s)-A'(s)}{1-A(s)}\rsz+\frac{R'(s)-d\cdot
  R(s)-C(s)}{1-A(s)}\,.  
\end{align*}
This is an inhomogeneous linear differential equation of order one for $G$. 
Restricted to the set of positive real numbers we know that the
solution is unique with any condition of type $G(s_0)=t_0$,  
and there is an explicit formula for it. Introducing the notations 
\begin{equation*}
L(s)=d-\frac{B(s)-A'(s)}{1-A(s)}\,,\quad
R^*(s)= -\frac{R'(s)-d\cdot R(s)-C(s)}{1-A(s)}\,, 
\end{equation*}
all solutions of the differential equation can be obtained in the
 form 
\begin{equation}\label{reklap2}
G(s)=\exp\biggl(-\int_s^{1} L(t)\,dt\biggr)\lsz C_0 +\int_s^{1}
R^*(t)\exp\biggl(\int_t^{1} L(u)\,du\biggr)dt\rsz 
\end{equation}
for $ s>0$, with an appropriate constant $C_0$.

From the results of Theorem \ref{rekt2} it follows that
$G(s)\rightarrow 0$ as $s$ goes to infinity on the real line. 
On the other hand, $L(t)\rightarrow d>0$ as $t\tv$, thus the first
exponential tends to infinity as $s\tv$. 
Hence there can exist at most one $C_0$ for which equation
\eqref{reklap2} is satisfied. Conditions (i5) and (i6) imply that
$R^*(t)\leq C_1/t$ for some $C_1$; in addition, $L(t)\leq C_2 d$ also
holds for $t>1/2$. Therefore we have 
\[
\int_s^{+\infty}R^*(t)\exp\biggl(\int_t^{s} L(u)\,du\biggr)dt\leq
\int_s^{+\infty} \frac{C_1}{t}\,\exp\bigl(C_2d(s-t)\bigr)dt
\leq\frac{C_3}{s}
\]
with some constant $C_3$, if $s>1/2$.

By this, setting 
\[
C_0=\int_1^{+\infty} R^*(t)\exp\lef
\int_t^1 L(u)\,du\ri\,dt,
\]
which is finite, in \eqref{reklap2} we get that
\begin{equation}\label{reklap3}
G(s)=\int_s^{+\infty} R^*(t)\exp\biggl(\int_t^{s}L(u)\,du\biggr)dt,
\end{equation}
for $s>0$, and this $G(s)\rightarrow 0$ as $s\tv$ on the real
line. Hence this is the Laplace transform of $g$ on the set of
positive numbers.  

Since $G(s)$ is holomorphic on the half-plane $\mathbb H$, it is the
unique extension of the solution given above. The right-hand 
side of \eqref{reklap3} is well-defined on $\mathbb H$, giving
a holomorphic function on $\mathbb H$, which extends $G$ from the set
of positive real numbers to $\mathbb H$. Thus we have that the Laplace
transform of $g$ is given by \eqref{reklap3} on the whole half-plane
$\mathbb H$.

Now we examine the behaviour of $G$ around zero. In what follows
$h_1,\,h_2,\,\dots$ will always denote functions that are holomorphic
in a neighbourhood of the origin. Let us start with $L$. 
Using the Taylor expansion of the exponential function we get
\begin{equation}\label{reklap5}
A(s)=1-s\int_0^{\infty} ta(t)\,dt+\frac{s^2}{2} \int_0^{\infty} t^2
a(t)\,dt+\ldots\,,\quad s\in\mathbb C, 
\end{equation}
which implies that 
\begin{equation*}
L(s)=d-\frac{B(s)-A'(s)}{1-A(s)}=\frac{B(0)-A'(0)}
{s\int_0^{\infty}ta(t)\,dt}+h_1(s)=
-\frac{\gamma+1}{s}+h_1(s).
\end{equation*}

Furthermore we have 
\begin{equation*}
\int_{s}^1 L(t)\,dt =(\gamma+1)\log s+h_2(s),\quad s\in\mathbb H. 
\end{equation*} 

Here we chose an arbitrary holomorphic branch of the logarithm on the
right half-plane $\mathbb H$.  
Once the logarithm is defined on $\mathbb H$, then $s^{\gamma+1}$ and
$s^{-(\gamma+1)}$ are also meaningful there. Thus we obtain that 
\begin{equation}\label{intL}
\exp\biggl(\int_t^{s}L(u)\,du\biggr)=\lef\frac ts\ri^{\gamma+1}\,
\exp\bigl(h_2(t)-h_2(s)\bigr),\quad s, t\in\mathbb H.
\end{equation}

One can similarly derive that $sR^*(s)$ is holomorphic in a
neighbourhood of $0$, hence
\begin{equation}\label{R*}
R^*(s)\exp\bigl(h_2(s)\bigr)=\frac{h_3(s)}{s}\,,\quad s\in\mathbb H.
\end{equation}

Finally, from equation \eqref{reklap3} we obtain that 
\begin{multline}\label{G}
G(s)=\int_s^{+\infty}R^*(t)\exp\biggl(\int_t^s L(u)\,du\biggr)\,
dt \\
=s^{-(\gamma+1)}\,\exp\bigl(-h_2(s)\bigr)\int_s^{+\infty}
t^{\gamma+1}R^*(t)\exp\bigl(h_2(t)\bigr)\,dt, \quad s\in\mathbb H. 
\end{multline}

Suppose first that $\gamma$ is not a negative integer, and consider only
positive values of $s$. Then, with a sufficiently small positive
$\veps$, by \eqref{R*} we have   
\begin{align*}
\int_s^{\infty}  t^{\gamma+1}R^*(t) \exp\bigl(h_2(t)\bigr)\,dt
&=C_4+\int_s^{\veps} t^{\gamma+1} R^*(t)\exp\bigl(h_2(t)\bigr)\,dt\\
&=C_4+\int_s^{\veps}t^\gamma h_3(t)\,dt\\
&=C_4+s^{\gamma+1}h_4(s)
\end{align*}
for $s\in(0,\veps)$, where $C_4$ is a constant. Hence
\[
G(s)=\exp\bigl(-h_2(s)\bigr)\left(C_4\,s^{-(\gamma+1)}+h_4(s)\right)
=h_5(s)+s^{-(\gamma+1)}h_6(s),
\]
from which the $k$th derivative of $G$ can be written in the following
form. 
\begin{equation*}
G^{(k)}(s)=h_7(s)+s^{-(\gamma+1+k)}h_8(s).
\end{equation*}
 
Choose a positive integer $k$ such that $0<\gamma+k+1$, then it
follows that
\begin{equation}\label{reklap6}
s^{\gamma+k+1}G^{(k)}(s)\rightarrow K
\end{equation}
as $s\to +0$, with some finite constant $K$.

Before going further, we prove a similar relation for $\gamma=-k$,
where $k$ is a positive integer. In this case we have
\begin{align*}
\int_s^{\infty}  t^{\gamma+1}R^*(t) \exp\bigl(h_2(t)\bigr)\,dt
&=C_4+\int_s^{\veps} t^{\gamma+1} R^*(t)\exp\bigl(h_2(t)\bigr)\,dt\\
&=C_4+\int_s^{\veps}t^\gamma h_3(t)\,dt\\
&=C_4+s^{\gamma+1}h_4(s)+C_5\log s,
\end{align*}
where $C_4$ and $C_5$ are constants, and $s\in (0,\veps)$. The term $\log s$ comes from the $(k-1)$st term of the expansion 
of the holomorphic function $h_3$. Then  
\begin{align*}
G(s)&=\exp\bigl(-h_2(s)\bigr)\left(h_9(s)+C_5\,s^{k-1}\log s\right)\\
&=h_{10}(s)+s^{k-1}h_{11}(s)\log s,
\end{align*}
consequently, $G^{(k-1)}(s)=h_{12}(s)+h_{13}(s)\log s$, and finally
\begin{equation*}
G^{(k)}(s)=\frac{h_{14}(s)}{s}+h_{15}(s)\log s, \quad s\in(0, \veps).
\end{equation*}
This implies that
\[
s\,G^{(k)}(s)\to K
\]
as $s\to +0$, with some finite constant $K$. Thus, \eqref{reklap6}
remains valid for negative integer values of $\gamma$.

Now we apply Karamata's Tauberian theorem (see e.g. 
\cite[Theorem XIII.5.2]{feller}, \cite[Theorem 1.7.1]{bingham}). We
will use the
following notation. Functions $v$ and $w$ are asymptotically equal to
each other, that is, $v\lef x\ri\sim w\lef x\ri$ as $x\rightarrow 0$
(or $\infty$) if $v/w$ tends to 1 as $x\rightarrow 0$ (or
$\infty$). $v\lef x\ri\sim 0\cdot w\lef x\ri$ means that $v/w$ tends
to 0 as $x\rightarrow 0$ (or $\infty$). The latter is the same as
$v=o\lef w\ri$.   

\begin{thmcite}
 Let $U$ be a non-decreasing right-continuous function on $\lef
 0,\infty\ri$ such that its Laplace transform
 $\omega(s)=\int_0^{\infty} e^{-sx}dU\lef x\ri$ exists for $s>0$. If
 $\ell$ is slowly varying at infinity, $0\leq \rho<\infty$, and $0\leq
 c<\infty$, then each of the relations  
\[
\omega(s)\sim cs^{-\rho}\ell\lef \frac{1}{s}\ri\quad\text{as } 
s\rightarrow +0
\] 
and
\[
U(t) \sim c\,\frac{1}{\Gamma(\rho+1)}\,t^{\rho}\ell(t)\quad 
\text{as }t\tv
\]
implies the other.
\end{thmcite}

We apply this theorem to $U\lef x\ri=\int_0^x g\lef u\ri u^k du$, for
which $\omega$ is  constant times $G^{(k)}$. From 
equation \eqref{reklap6} we get that
\begin{equation}\label{reklap6.2}
\int_0^x g\lef u\ri u^k du\sim A_kx^{\gamma+k+1} 
\end{equation}
as $x\tv$, for some  $A_k\geq 0$. Note that the constant $A_k$ depends on $k$. 

In order to finish the proof of Theorem \ref{rekt3} we need another
Tauberian type theorem, giving the asymptotics of $g(u)u^k$  from the
asymptotics of its integral function. We will use the monotonicity of
$g$ at this point.  

We say that a function $f$ is slowly oscillating if  for any $\veps>0$
there exists $\delta>0$ such that $f\lef u\ri<f\lef x\ri\lef
1+\veps\ri$ holds for all $x<u<x\lef 1+\delta\ri$ (see
e.g. \cite[Section 6.2]{hardy}, \cite[Section 1.7.6,]{bingham},
\cite[Section 17]{korevaar}). Using that
$g$ is a non-increasing, nonnegative function, it is easy to see that
$g\lef x\ri x^k$ is slowly oscillating provided $k$ is a positive
integer. Indeed, for all $x<u<x\lef 1+\delta\ri$ we have  
\[
g\lef u\ri u^k \leq g\lef x\ri x^k \lef 1+\delta\ri^k.
\]
Hence given $\veps>0$, any $\delta>0$ such that $\lef
1+\delta\ri^k<1+\veps$ would satisfy the condition.

Slow oscillation is generally  a sufficient condition of Tauberian type
theorems. For example, Theorem 17.2. of \cite{korevaar} states the
following.  

\begin{thmcite1}
Let $f$ be defined on an interval $\lef a,\infty\ri$, and suppose that
$f(x)\sim Ax^{\al}$ as $x\tv$ with some real numbers $\al,\,A$. 
If $f$ is $m$ times differentiable and 
\[
\liminf\frac{f^{(m)}(y)-f^{(m)}(x)}{x^{\al-m}}\geq 0
\]
as $x\tv$ and $1<y/x\to 1$, then 
\[
f^{(j)}(x)\sim A\al(\al-1)\dots(\al-j+1)\,x^{\al-j}
\]
as $x\tv$, for $1\leq j\leq m$.
\end{thmcite1}

Based on equation \eqref{reklap6.2}, we can apply this theorem with 
\[
f(x)=-\int_0^x g(u)u^kdu,\quad m=j=1,\quad\text{and}
\quad\al=\gamma+k+1. 
\]
Then we get that $g\lef x\ri x^{-\gamma}$ is convergent as $x\tv $.

From Theorem \ref{rekt2} it follows that the limit of $g\lef x\ri
x^{-\gamma}$ is positive and finite, which is just our Theorem
\ref{rekt3}. \qed

\begin{remark}
Since the Laplace transform method usually gives only local 
results in the discrete case, and we needed global results there, 
it is reasonable to use classic renewal techniques. 
On the other hand, those methods rely on convolution 
in the continuous case, which was not useful for our integral equation.
\end{remark}

\subsubsection*{Acknowledgement} Authors are indebted to G\'abor Hal\'asz for his invaluable help with Laplace transforms.

\end{document}